\documentclass[conference,a4paper]{IEEEtran}
\IEEEoverridecommandlockouts
\usepackage{lipsum}
\usepackage{graphicx}
\usepackage{epstopdf}
\usepackage{algorithmic}
\usepackage{amsmath} 
\usepackage{amsfonts,mathrsfs}
\allowdisplaybreaks 
\usepackage{amssymb}
\usepackage{color}
\usepackage{graphicx}
\usepackage{epsfig}
\usepackage{subfigure}
\usepackage{enumerate}
\usepackage{algorithm}
\usepackage{empheq}
\usepackage{bm}
\usepackage{accents}
\usepackage{balance}
\usepackage{marginnote}
\usepackage{setspace}
\usepackage{titlesec}
\usepackage{url}
\usepackage{soul}
\def\BibTeX{{\rm B\kern-.05em{\sc i\kern-.025em b}\kern-.08em
		T\kern-.1667em\lower.7ex\hbox{E}\kern-.125emX}}

\newcommand{\mc}{\mathcal}
\newcommand{\bb}{\mathbb}

\DeclareMathAlphabet{\mathbbmsl}{U}{bbm}{m}{sl}

\newcommand{\col}{\operatorname{col}}

\newcommand{\1}{\mathbf{1}}

\newcommand{\Rmnum}[1]{\expandafter\@slowromancap\romannumeral #1@}

\begin{document}
	
	\title{Operating envelopes for the grid-constrained use of distributed flexibility in balancing markets
		\thanks{This work is supported by the Belgian FPS economy through the Energy Transition Funds project ALEXANDER.}\vspace{-0.4cm}
	}
	
	\author{\IEEEauthorblockN{Abhimanyu Kaushal, Wicak Ananduta, Luciana Marques, Tom Cuypers, and Anibal Sanjab}
		\IEEEauthorblockA{
			\textit{Flemish Institute for Technological Research (VITO) and EnergyVille}, Genk, Belgium \\
			\{abhimanyu.kaushal, wicak.ananduta, luciana.marques, tom.cuypers, anibal.sanjab\}@vito.be\vspace{-0cm}}
	}
	\maketitle
	\begin{abstract}
		The increasing share of distributed energy sources enhances the participation potential of distributed flexibility in the provision of system services. However, this participation can endanger the grid-safety of the distribution networks (DNs) from which this flexibility originates.
		In this paper, the use of operating envelopes (OE) to enable the grid-safe procurement of distributed flexibility in centralized balancing markets is proposed. 
		Two classes of approaches for calculating OEs (one-step and two-step methods) are compared in terms of the level of distribution grid safety they can provide, the impact they can have on the market efficiency, and the volume of discarded flexibility they can yield. A case study considering different system scenarios, based on Monte Carlo simulations, highlights a trade-off between the market efficiency, DN flexibility resource utilization, and the grid safety delivered by the different OE methods. The results showcase that the use of the two-step OE approach results in a more grid-secure albeit less-efficient use of distributed flexibility.
		
	
	\end{abstract}
	
\section{Introduction}
The {large-scale integration} of {variable} renewable energy resources (RES) {and the expected increase in load levels driven by advanced electrification} increase the system's flexibility requirements~\cite{JRC130519}. 
The flexibility availability from prosumers/resources in the distribution network (DN) would then be crucial for the economic and secure operation of the power system 
~\cite{EUDSOcommissionregulation}.
This enhanced flexibility requirement along with the technological evolution presents the prosumers connected at the DN level with an economic opportunity to adjust their energy consumption for providing flexibility to transmission or distribution system operators (TSO/DSO). Indeed, the work in~\cite{marques2023grid} have presented different market-clearing models for flexibility procurement by TSOs of DN resources. 
{The use of DN flexibility assets in transmission-level services requires the inclusion of DN constraints in the transmission-level market clearing problem to ensure that the DN operational limits are respected when activating distribution-level flexibility.} 
However, this inclusion of {full} DN constraints faces hurdles in terms of communication efficiency ({due to the increased burden of continuously updating and communicating the network models}) and privacy preservation {(in terms of the need by the DSO for sharing DN constraints externally).}
To {impose distribution-grid limitations without the need for a full network representation, the DSO can fix} import and export limits at connection points~\cite{petrou2021ensuring}. 
However, as the potential of flexibility from DN resources increases, these fixed limits {can over} restrict the power exchange from DN resources, thus, unnecessarily harming the profitability of prosumers~\cite{petrou2021ensuring} which, in turn, can hinder their participation in flexibility markets.

The use of operating envelopes (OEs) as dynamic limits for flexibility from DN resources has been presented as a potential solution to utilize the full potential of the DN assets while accounting for grid safety in~\cite{petrou2021ensuring, liu2021grid, liu2022using, 8810638}. 
The OEs can represent the dynamic flexibility limits for individual resources~\cite{petrou2021ensuring, liu2021grid, liu2022using} or for a {complete} DN (limit at the TSO/DSO interface)~\cite{8810638}. 
The flexibility limits resulting from the OE approaches (per resource/DN) can also be used in TSO-level flexibility market formulations, thus replacing the need for communicating DSO network models. For example, the use of OEs for peer-to-peer trading is presented in~\cite{azim2024dynamic, hoque2024dynamic}. 
However, {the utilization of OE methods for using DN resource flexibility for TSO-level flexibility market clearing is not explored yet and at the same time a comparison of using different OE calculation methods is also missing.} 




This paper introduces and analyzes the suitability of using {OEs} for the participation of distributed flexibility in TSO-level balancing markets. The analysis focuses on their resulting impact on market efficiency (cost of flexibility procurement), grid security (number of constraint violations), and conservativeness (limitation on the resource flexibility participation).  
%
%
{Two} different OE calculation approaches (one and two-step methods) are presented and the influence of flexibility resource prices and their quantities on the calculated OE volume is examined.
%
%
To analyze the performance of the different OEs, a case study is performed capturing differing settings in terms of availability of generation at the DN level and the number of available DN flexibility bids. A Monte Carlo based approach is used to capture a wide range of flexibility volume and price scenarios. 
\subsection*{Notation} 
The set of real numbers is denoted by $\bb R$. The operator $\col(\cdot)$ concatenates its arguments column-wise. The vector whose all elements are 1 is denoted by $\1_n \in \bb R^n$. A set is indicated by a calligraphic uppercase letter and we add a subscript if it is a subset.
The set of distribution networks is denoted by $\mc D:=\{d_1, d_2, \dots, d_M\}$, for $M>0$. For the transmission network (TN), subscript $T$ indicates variables, parameters, and sets. So, the set of all networks is denoted by $m \in \mc N:=\mc D \cup \{T\}$.  For each network $m \in \mc N$, $\mc B_m$ and $\mc L_m$ denotes the set of busses and lines, respectively. 
\section{Central flexibility market model} \label{sec:bid_forwarding}
\vspace{-3pt}
We consider a balancing market, e.g., frequency containment or manual frequency restoration reserve, that allows participation of DN flexibility service providers. To make sure that any flexibility activation do not cause congestion in the transmission system, the network constraints of transmission system are also considered.
To show the formulation of the market clearing problem, let us first denote by $\mc R:= \{r_1, r_2, \dots r_N\}$ the set of all flexibility resources. For each $n \in \mc R$, the flexibility quantity (in terms of active power) is denoted by $p_n \in \bb R$. 
For each network, $m \in \mc N$, the concatenation of the network flexibility quantity is denoted by $\bm p_m := \col((p_n)_{n \in {\mc R_m}})$. The interface flow between a distribution system $m \in \mc D$ and the transmission bus to which $m$ is connected is denoted by $z_m$ and we define $\bm z := \col((z_m)_{m \in \mc D})$. The collection of power injection ($P_i^{\mathrm{in}}$) at all transmission busses is denoted by $\bm \varphi_T := \col((P_i^{\mathrm{in}})_{i \in \mc B_T})$.
We can then cast the market clearing problem as \vspace{-0.2cm}

{\small
\begin{subequations}
	\begin{align}
		& \min_{(\bm p_T, \bm \varphi_T), (\bm p_m, z_m)_{ \forall {m \in \mc D}}} \ \  \sum_{n \in \mc R} c_n p_n  \label{eq:objective_balancing_market}\\
		\operatorname{s.t.} \ \ & 
		\underline{p}_n \leq p_n \leq \overline p_n, \ \  \forall n \in \mc R, \label{eq:resource_constraints}\\
		& P_i^{\mathrm{in}} = \sum_{n \in \mc R_i} p_n + e_i, \quad \forall i \in \mc B_T, \label{eq:nodal_balance_transmission}\\
		&  - \overline{\bm \varphi}_T^{\mathrm{f}} \leq C (\bm \varphi_T - D \bm z ) \leq \overline{\bm \varphi}_T^{\mathrm{f}}, \label{eq:PTDF_transmission_constraints} \\
		& z_m + \sum_{n \in \mc R_m} \bm p_n = {z_m^0}, \quad \ \forall m \in \mc D, \label{eq:aggregated_distribution_balance_constraints} \\
		& -\overline{z}_m \leq z_m \leq \overline z_m,\quad \ \forall m \in \mc D. \label{eq:limit_interface_flows}
	\end{align}%
	\label{eq:balancing_market_problem}%
\end{subequations}
}%
The objective of Problem \eqref{eq:balancing_market_problem} is to minimize the 
flexibility procurement cost defined in \eqref{eq:objective_balancing_market} where $c_n > 0$ denotes the price of individual bids.  
The constraints in \eqref{eq:resource_constraints} define the range of flexibility that can be offered by each resource with $\underline{p}_n$ and $\overline{p}_n$ as the minimum and maximum values. 
We assume that each resource $n \in \mc R$ can only provide an upward flexibility, i.e., $\underline{p}_n=0$ and $\overline{p}_n > 0$, or a downward flexibility, i.e., $\underline{p}_n<0$ and $\overline{p}_n = 0$. We denote by $\mc R^{\mathrm u}$ and $\mc R^{\mathrm d}$ the set of upward resources and that of downward resources, respectively. The constraints in \eqref{eq:nodal_balance_transmission} represent the nodal balance and \eqref{eq:PTDF_transmission_constraints} defines the power flow constraints of the TN with $\overline{\bm \varphi}_T^{\mathrm{f}} \in \bb R_{>0}^{|\mc L_T|}$ being the maximum limit. Note that we consider the PTDF model where the power flows are linearly proportional to the power injections with $C$ being the PTDF matrix and $D$ being a selection matrix consisting of zeros and ones, i.e., the element $[D_0]_{i,m} = 1$ if distribution system $m$ is connected to transmission bus $i$, otherwise $[D_0]_{i,m} = 0$. The balancing equations of all DNs is shown in~\eqref{eq:aggregated_distribution_balance_constraints}, where $z_m^0$ is the initial interface flow (the net of anticipated base injections and loads in the entire network).  Finally, \eqref{eq:limit_interface_flows} defines the the interface flow constraints with $\overline{z}_m$ being the maximum limit. 

As DN operational constraints are not considered in Problem \eqref{eq:balancing_market_problem}, unfortunately the market operator cannot guarantee that if any distribution-level flexibility bids are activated, they will not cause network issues such as congestion since they do not abide by the network constraints. The existing literature on such market formulations, e.g., \cite{marques2023grid, sanjab2023joint, marques2023flexibility}, indeed emphasizes the inclusion of a DN model. Nevertheless, in practice, this might not always be possible due to the sensitivity and/or unavailability of data that must be obtained by the market operator from the DSOs. 
\section{Operating envelope integration} \label{sec:oe_calculation}
To overcome the aforementioned grid security issue, this section
presents the integration of the OE approach to the market formulation and the methodologies to calculate OEs.
\vspace{-3pt}
\subsection{Integration}
An operating envelope, denoted by $\mc E_n$, is defined as a feasible range in which a resource may operate, i.e.,

{
\begin{equation}
	\mc E_n := \{p \in \bb R \mid \underline{\varepsilon}_n \leq p \leq \overline{\varepsilon}_n\}, 
	\end{equation}}
	
	\noindent which are parameterized by the lower $\underline{\varepsilon}_n$ and upper $\overline{\varepsilon}_n$ limits. We further denote by $\underline{\bm \varepsilon}_m = \col((\underline{\varepsilon}_n)_{n \in \mc R_m})$ and $\overline{\bm \varepsilon}_m = \col((\overline{\varepsilon}_n)_{n \in \mc R_m})$ the concatenation of these resource limits for each distribution system $m \in \mc D$. 
	Supposing that each distribution-level resource $n \in \mc R_m$, for each $m \in \mc D$, has an OE $\mc E_n$, the market clearing problem in \eqref{eq:balancing_market_problem} turns into\vspace{-0.2cm}
	
	{
\begin{subequations}
	\begin{align}
		& \min_{(\bm p_T, \bm \varphi_T), (\bm p_m, z_m)_{ \forall {m \in \mc D}}} \ \  \sum_{n \in \mc R} c_n p_n  \label{eq:objective_balancing_market2}\\
		\operatorname{s.t.} \ \ & 
		\text{\eqref{eq:resource_constraints}--\eqref{eq:limit_interface_flows} and} \notag \\
		& p_n \in \mc E_n,\ \forall n \in \mc R_m,\ \forall m \in \mc D, \label{eq:oe_constraints}
	\end{align}%
	\label{eq:balancing_market_problem_OE}%
\end{subequations}\vspace{-0.4cm}
}%

\noindent where  we enforce the DN resources to be in their OEs by \eqref{eq:oe_constraints}. 
These OEs in the market clearing process can be considered as a prequalification step with the objective of ensuring only safe (amount of) bids can be submitted to the market.
\subsection{Calculation methods}

Next, we review and generalize the existing calculation methods for computing an OE.
These methods rely on solving at least an optimization problem that takes into account the network constraints,
which couple $\bm p_m$ with physical variables of the network (concatenated as $\bm \varphi_m$). Specifically, in this work, we consider the linearized Branch Flow model \cite[Eqs. (3)--(7)]{sanjab2021linear}, where $\bm \varphi_m := \col((P_i^{\mathrm{in}}, Q_i^{\mathrm{in}}, v_i)_{i \in \mc B_m},(P_{(i,j)}^{\mathrm{f}}, Q_{(i,j)}^{\mathrm{f}})_{(i,j) \in \mc E}, z_m^{\mathrm{re}})$, with $P_i^{\mathrm{in}}, Q_i^{\mathrm{in}}, v_i$ denoting the nodal active power injection, reactive power injection, and voltage, $P_{(i,j)}^{\mathrm{f}}, Q_{(i,j)}^{\mathrm{f}}$ denoting the active and reactive power flows, and $z_m^{\mathrm{re}}$ denoting the reactive interface flow. We assume that each distribution system has a radial structure and we denote by $\pi(i)$ the parent node of bus $i \in \mc B_m$ and by $b_{0,m}$ the root node. The model is defined by the following constraints:

{\small
\begin{align}
	&P_i^{\mathrm{in}} = 
	\sum_{n \in \mc R_i} p_n + e_i, \quad \forall i \in \mc B_m.
	\label{eq:injected_power}\\
	&Q_i^{\mathrm{in}} = 
	\sum_{n \in \mc R_i} \alpha_n p_n + e_i^{\mathrm{re}}, \quad \forall i \in \mc B_m, \label{eq:injected_power_reactive} \\
	&P_i^{\mathrm{in}} = 
	\begin{cases}
		\sum_{j \in \mc K_i} P_{(i,j)}^{\mathrm{f}} - z_m, & i=b_{0,m}, \\
		\sum_{j \in \mc K_i} P_{(i,j)}^{\mathrm{f}} - P_{(\pi(i), i)}^{\mathrm{f}} , & i \in \mc B_m \setminus\{b_{0,m}\},
	\end{cases}
	\label{eq:lindistflow_1} \\
	&	Q_i^{\mathrm{in}} = 
	\begin{cases}
		\sum_{j \in \mc K_i} Q_{(i,j)}^{\mathrm{f}} - z_m^{\mathrm{re}}, & i =b_{0,m}, \\
		\sum_{j \in \mc K_i} Q_{(i,j)}^{\mathrm{f}} - Q_{(\pi(i), i)}^{\mathrm{f}} , & i \in \mc B_m \setminus\{b_{0,m}\},
	\end{cases}
	\label{eq:lindistflow_2} \\
	&	v_i = 
	\begin{cases}
		v_m^0, \qquad  i =b_{0,m}, \\
		v_{\pi(i)} \hspace{-2pt}-\hspace{-2pt} 2 (R_{(\pi(i),i)} P_{(\pi(i),i)}^{\mathrm f} \hspace{-2pt}+\hspace{-2pt} X_{(\pi(i),i)} Q_{(\pi(i),i)}^{\mathrm f}), \\ \qquad \qquad \forall i \in \mc B_m\setminus\{b_{0,m}\}, 
	\end{cases}	
	\label{eq:v_lindistflow} \\
	& E^{P} P_{(i,j)}^{\mathrm{f}} + E^{Q} Q_{(i,j)}^{\mathrm{f}} \leq E^{S} S_{(i,j)}^{\mathrm{f}} \label{eq:flow_bound} \\
	& \underline{v}_i \leq v_i \leq \overline{v}_i, \quad \forall i \in \mc B_m, \label{eq:voltage_bound_lindistflow} \\
	& \underline{z}_m^{\mathrm{re}} \leq z_m^{\mathrm{re}} \leq \overline{z}_m^{\mathrm{re}}, \label{eq:reactive_interface_lindistflow_bound}
\end{align}
}

\noindent where \eqref{eq:injected_power} is the nodal active power balance equation with $e_i$ being the net of anticipated base active power injection and load (note that $\sum_{i \in \mc B_m} e_i = z_m^0$); \eqref{eq:injected_power_reactive} is the nodal reactive power balance equation assuming that the reactive power of each distribution-level resource is linearly proportional to the active power with $\alpha_n>0$ being the linear proportion and  $e_i^{\mathrm{re}}$ being the net of anticipated base reactive power injection and load; \eqref{eq:lindistflow_1}--\eqref{eq:lindistflow_2} are the power flow equations in the network; \eqref{eq:v_lindistflow} describes the voltage equation where $v_m^0$ denotes the operating voltage of distribution system $m \in \mc D$; 
\eqref{eq:flow_bound} is a linear inner approximation of the quadratic power flow bound $S_{(i,j)}^{\mathrm f}$ parameterized by the matrices $E^{P}$, $E^{Q}$, and $E^{S}$;
\eqref{eq:voltage_bound_lindistflow} bounds the bus voltages; and \eqref{eq:reactive_interface_lindistflow_bound} bounds the reactive power interface flow. Therefore, we can define the network constraint set as:\vspace{-0.2cm}

{
\begin{equation}
	\begin{aligned}
		\mc C_m := \left\{(\bm p_m, \bm \varphi_m, z_m) \mid \text{\eqref{eq:injected_power}--\eqref{eq:reactive_interface_lindistflow_bound} hold}  \right\}.
	\end{aligned}
	\label{eq:lindistflow_constraint_set}
\end{equation}
}

We study two different OE calculation methods presented in \cite{liu2022using} and \cite{petrou2021ensuring}, and shown in Algorithms~\ref{alg:oe_liu2022} and \ref{alg:oe_petrou2021}, respectively. We note that these algorithms are performed individually by each DSO. Algorithm~\ref{alg:oe_liu2022} is a two-step approach, i.e., the calculation of the limits of the upward and downward resources are done separately by solving Problems \eqref{eq:oe_problem_liu2022_max} and \eqref{eq:oe_problem_liu2022_min}, respectively. Notice that in Problems \eqref{eq:oe_problem_liu2022_max} and \eqref{eq:oe_problem_liu2022_min}, the network constraint $\mc C_m$ as defined in \eqref{eq:lindistflow_constraint_set} is included (e.g., see \eqref{eq:dn_constraint}). Furthermore,  in Problem \eqref{eq:oe_problem_liu2022_max}, all downward resources are set to 0 by \eqref{eq:oe_problem_liu2022_max_resource_constraint} and, similarly, in Problem \eqref{eq:oe_problem_liu2022_min}, all upward resources are set to 0 by \eqref{eq:oe_problem_liu2022_min_resource_constraint}. On the other hand, Algorithm~\ref{alg:oe_petrou2021} calculates the limits of all resources simultaneously by solving Problem \eqref{eq:oe_problem_quadratic_obj}. These two approaches in general consider the same constraints, i.e., the network constraints and the operational constraints of the resources. On the other hand, the objective functions used in these methods are different. Namely, in \cite{liu2022using}, the objective of each optimization is to maximize the range of the OEs by considering a linear function, i.e., maximizing $p_n$ for computing $\overline{\varepsilon}_n$ (of the upward resources), as in \eqref{eq:oe_problem_liu2022_max_obj}, and minimizing $p_n$ for computing $\underline{\varepsilon}_n$ (of the downward resources), as in \eqref{eq:oe_problem_liu2022_min_obj}. Meanwhile, in \cite{petrou2021ensuring}, a quadratic function is considered since the objective is to find the closest operating point to the limit, i.e., $\overline{p}_n$ for upward resources and $\underline{p}_n$ for downward resources.   

We generalize the objective functions considered in these two methods by introducing weights, denoted by $w_n$, whose purpose is to allow for having preferences/priorities for certain resources, depending on their prices, availability, etc. The original formulations in \cite{liu2022using} and \cite{petrou2021ensuring} use equal weights, i.e., $w_n = 1$, for all $n \in \mc R_m$ and all $m \in \mc D$. However, two other different weights can be considered, namely price-based, quantity-based. Under the price-based weight, we set\vspace{-0.2cm}

{\small
\begin{equation}
	w_n = \begin{cases}
		{c_n}/{(\max_{n \in \mc R_m} c_n)}, \ \ \forall n \in \mc R_m^{\mathrm d}, \\
		{(\max_{n \in \mc R_m} c_n)}/{c_n}, \ \ \forall n \in \mc R_m^{\mathrm u},
	\end{cases}  \forall m \in \mc D.
	\label{eq:price_based_weight}
\end{equation}
}

Here, the prices are normalized by the most expensive resources. This rule aims at giving the highest priorities (i.e., the highest $\underline{\varepsilon}_n$ and $\overline{\varepsilon}_n$) to the cheapest resources, aligning with the market clearing formulation in its objective to be most efficient. In this case, DSOs need to know the prices of the resources.
On the other hand, the quantity-based weights,

{\small
\begin{equation}
	w_n = \max\{\overline{p}_n, -\underline{p}_n \}, \ \ \forall n \in \mc R_m, \forall m \in \mc D,
	\label{eq:quantity_based_weight}
\end{equation}
}

\noindent prioritize resources with the largest quantity. This weight rule does not need any other additional information from the resource owners as the parameters used are the same as those used in the constraint formulation. 
\begin{algorithm}[!t]
\small
\caption{Two-step operating envelope calculation}
\label{alg:oe_liu2022}
\begin{enumerate}
	\item Compute $(\bm p_m^*, \bm \varphi_m^*, z_m^*)$ as a solution to 	
	\begin{subequations}
		\begin{align}
			& \underset{\bm p_m, \bm \varphi_m, z_m}{\max} \ \  \sum_{n \in \mc R_m} {w_n}p_n \label{eq:oe_problem_liu2022_max_obj}\\
			& \operatorname{s.t.} \ \  
			0 \leq p_n \leq \overline{p}_n, \qquad \quad \ \  \forall n \in \mc R_m, \label{eq:oe_problem_liu2022_max_resource_constraint} \\
			&\qquad  \ (\bm p_m, \bm \varphi_m, z_m) \in \mc C_m,  \label{eq:dn_constraint}\\
			& \qquad \   -\overline{z}_m \leq z_m \leq \overline z_m, \notag
		\end{align}%
		\label{eq:oe_problem_liu2022_max}%
	\end{subequations}
	and set $\overline{\bm \varepsilon}_m= \bm p_m^*$.
	\item Recompute $(\bm p_m^*, \bm \varphi_m^*, z_m^*)$ as a solution to 	
	\begin{subequations}
		\begin{align}
			& \underset{\bm p_m, \bm \varphi_m, z_m}{\min} \ \  \sum_{n \in \mc R_m} {w_n}p_n \label{eq:oe_problem_liu2022_min_obj}\\
			& \operatorname{s.t.} \ \  
			\underline{p}_n \leq p_n \leq 0, \qquad \quad \ \  \forall n \in \mc R_m, \label{eq:oe_problem_liu2022_min_resource_constraint}  \\
			&\qquad  \ (\bm p_m, \bm \varphi_m, z_m) \in \mc C_m,  \notag\\
			& \qquad \   -\overline{z}_m \leq z_m \leq \overline z_m, \notag
		\end{align}%
		\label{eq:oe_problem_liu2022_min}%
	\end{subequations}
	and set $\underline{\bm \varepsilon}_m= \bm p_m^*$.
\end{enumerate}
Outputs: $(\underline{\varepsilon}_n, \overline{\varepsilon}_n)$, for all $ n\in\mc R_m$.
\end{algorithm} 
\begin{algorithm}[!t]
\small
\caption{One-step operating envelope calculation}
\label{alg:oe_petrou2021}
\begin{enumerate}
	\item Set $\underline \varepsilon_n= \underline{p}_n = 0$, for all $n \in \mc R_m^{\mathrm u}$, and $\overline \varepsilon_n= \overline{p}_n = 0$, for all $n \in \mc R_m^{\mathrm d}$, and 
	compute $(\bm p_m^*, \bm \varphi_m^*, z_m^*)$ as a solution to 	
	\begin{equation}
		\begin{aligned}
			& \underset{\bm p_m, \bm \varphi_m, z_m}{\min} \ \  \sum_{n \in \mc R_m^{\mathrm u}} \hspace{-5pt}w_n (p_n - \overline{p}_n)^2 + \sum_{n \in \mc R_m^{\mathrm d}}\hspace{-5pt}w_n (p_n - \underline{p}_n)^2 \\
			& \operatorname{s.t.} \ \  
			\underline{p}_n \leq p_n \leq \overline p_n,  \ \  \forall n \in \mc R_m, \\
			&\qquad  \ (\bm p_m, \bm \varphi_m, z_m) \in \mc C_m,  \\
			& \qquad \   -\overline{z}_m \leq z_m \leq \overline z_m,
	\end{aligned}%
	\label{eq:oe_problem_quadratic_obj}%
\end{equation}
to set $\overline \varepsilon_n= p_n^*$, for all $n \in \mc R_m^{\mathrm u}$, and $\underline \varepsilon_n= p_n^*$, for all $n \in \mc R_m^{\mathrm d}$.
\end{enumerate}
Outputs: $(\underline{\varepsilon}_n, \overline{\varepsilon}_n)$, for all $ n\in\mc R_m$.
\end{algorithm}


\section{Case study} 
\label{sec:test_system_results}\vspace{-0.1cm}

{We systematically compare the performance of the two OE methods when they are integrated into the flexibility market clearing via Monte Carlo-based numerical simulations. 
We consider a network consisting of the IEEE 14-bus transmission system connected with the Matpower 69-bus and 141-bus distribution systems.
In each case instance, we assume that there is an imbalance in the transmission system, which is resolved by a balancing market. The imbalances in the transmission systems are obtained by randomly varying the loads (and the generation for the transmission system) while ensuring that there is no anticipated congestion in the DNs before the flexibility market is cleared. The flexibility resources in the TN and DNs are then generated by randomly determining their locations, maximum quantities ($\underline p_n, \overline p_n$), and prices ($c_n$). To be able to observe the performance of the OE methods, the flexibility prices of the DN resources are considered to be cheaper than those in the TN.
For the resources of DN $m \in \mc D$, the prices are randomly set as $c_n \sim U(35,55)$, for all $n \in \mc R_m^{\mathrm u}$, and $c_n \sim U(14,34)$, for all $n \in \mc R_m^{\mathrm d}$, while the price of the transmission-level resources are set $c_n \sim U(65,75)$, for all $n \in \mc R_T^{\mathrm u}$, and $c_n \sim U(1,11)$, for all $n \in \mc R_T^{\mathrm d}.$\footnote{$U(a,b)$ denotes a uniform distribution where $a$ and $b$ are the minimum and maximum values.}


In each Monte Carlo instance, not only we solve Problem \eqref{eq:balancing_market_problem_OE} where the OEs are calculated by Algorithms~\ref{alg:oe_liu2022} (one-step) and~\ref{alg:oe_petrou2021} (two-step) (the OE-based market clearing problems) but we also solve Problem \eqref{eq:balancing_market_problem}, where no DN constraints are included (labeled as the \emph{no-DN} market problem), as well as the \emph{full-DN} market clearing problem which includes all DN constraints, i.e., Problem \eqref{eq:balancing_market_problem_OE} where the constraint in \eqref{eq:oe_constraints} is substituted by \eqref{eq:dn_constraint}, for all $m \in \mc D$. The latter two are used as benchmarks. Furthermore, for each OE calculation method, we consider three weight rules, namely equal weights ($w_n=1$), price-based weights as in \eqref{eq:price_based_weight}, and quantity-based weights as in \eqref{eq:quantity_based_weight}. After clearing these markets, the cleared flexibility is activated and we check the states of the DNs by solving the power flow equations \eqref{eq:injected_power}--\eqref{eq:lindistflow_2}.

For the quantitative comparison, we consider the following key performance metrics:
\begin{enumerate}[i.]
\item \textit{Grid safety of the cleared flexibility}, measured by the total number of nodal voltage and branch flow values that violate their bounds,
\item \textit{Market (in)efficiency}, in terms of the procurement cost. We normalize it with the idealized market cost, i.e.,

{\small
	\begin{equation}
		\eta = \frac{\sum_{n \in \mc R} (c_n p_n - c_n p_n^{\circ})}{|\sum_{n \in \mc R} c_n p_n^{\circ}|} \times 100\%, \label{eq:inefficiency}
	\end{equation}
}

\noindent where $\eta$ denotes the inefficiency while $p_n^{\circ}$, for all $n \in \mc R$, denote a solution to the idealized market clearing problem.
\item \textit{Unqualified flexibility}, defined by the total amount of flexibility quantity that is excluded from the market as a result of the OE constraints, i.e.,\vspace{-0.2cm}

{\small
	\begin{align*}
		\delta^{\mathrm d} = \frac{\sum_{n \in \cup_{m \in \mc D} \mc R^{\mathrm u
			}} |\underline{p}_n - \underline \varepsilon_n|}{\sum_{n \in \cup_{m \in \mc D} \mc R^{\mathrm u
			}} |\underline{p}_n|} \times 100\%,\\
		\delta^{\mathrm u} = \frac{\sum_{n \in \cup_{m \in \mc D}\mc R_m^{\mathrm u}} (\overline{p}_n - \overline \varepsilon_n)}{\sum_{n \in \cup_{m \in \mc D} \mc R_m^{\mathrm u}} \overline{p}_n} \times 100\%, 
	\end{align*}
}
\end{enumerate}
}

\noindent where $\delta^{\mathrm d}$ and $\delta^{\mathrm u}$ denote the unqualified downward and upward flexibility, respectively.

We create two sets of cases, specified as follows:
\begin{enumerate}
\item Case set 1: The TSO has downward balancing needs and the DNs have shiftable loads as flexibility resources.
\item Case set 2: The TSO has upward balancing needs and the DNs have shiftable loads and distributed generation.
\end{enumerate}

\begin{figure}
\centering
\includegraphics[width=0.95\columnwidth]{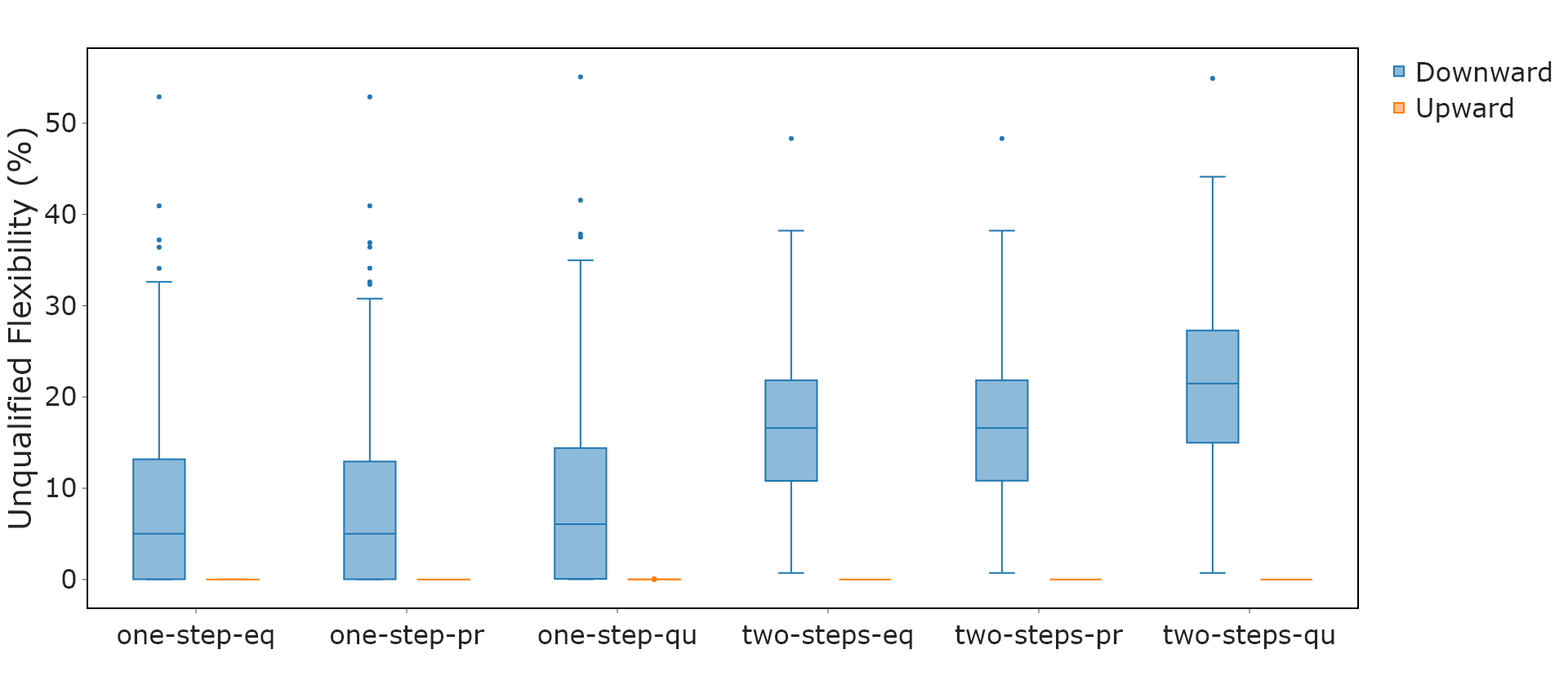}\vspace{-0.3cm}
\caption{Comparison of unqualified flexibility for case set 1.}
\label{fig:unqualified_flex_case1}\vspace{-0.4cm}
\end{figure}

The number of flexibility resources in case set 2 is on average $160\%$ more than that in case set 1. One can consider case set 1 as an illustration of the current situation while case set 2 represents scenarios in the future where the DNs have a high flexibility potential. In addition, for case set 1 we obtain $\sim$300 instances where a solution to the relaxed market clearing problem (\emph{no-DN}) is not grid-safe, i.e., there exist grid constraint violations\footnote{The instances where a solution to the relaxed market problem is grid-safe have been discarded as they are not interesting because the OEs further restrict the resources and thus solutions to the OE-based markets are also grid-safe.}. For case set 2, we obtain $\sim$1,600 instances. More instances are obtained for case 2 because this case has more distribute resources available, which can lead to more grid violations in the (\emph{no-DN}) scenario.
\begin{figure}
\centering
\includegraphics[width=\columnwidth]{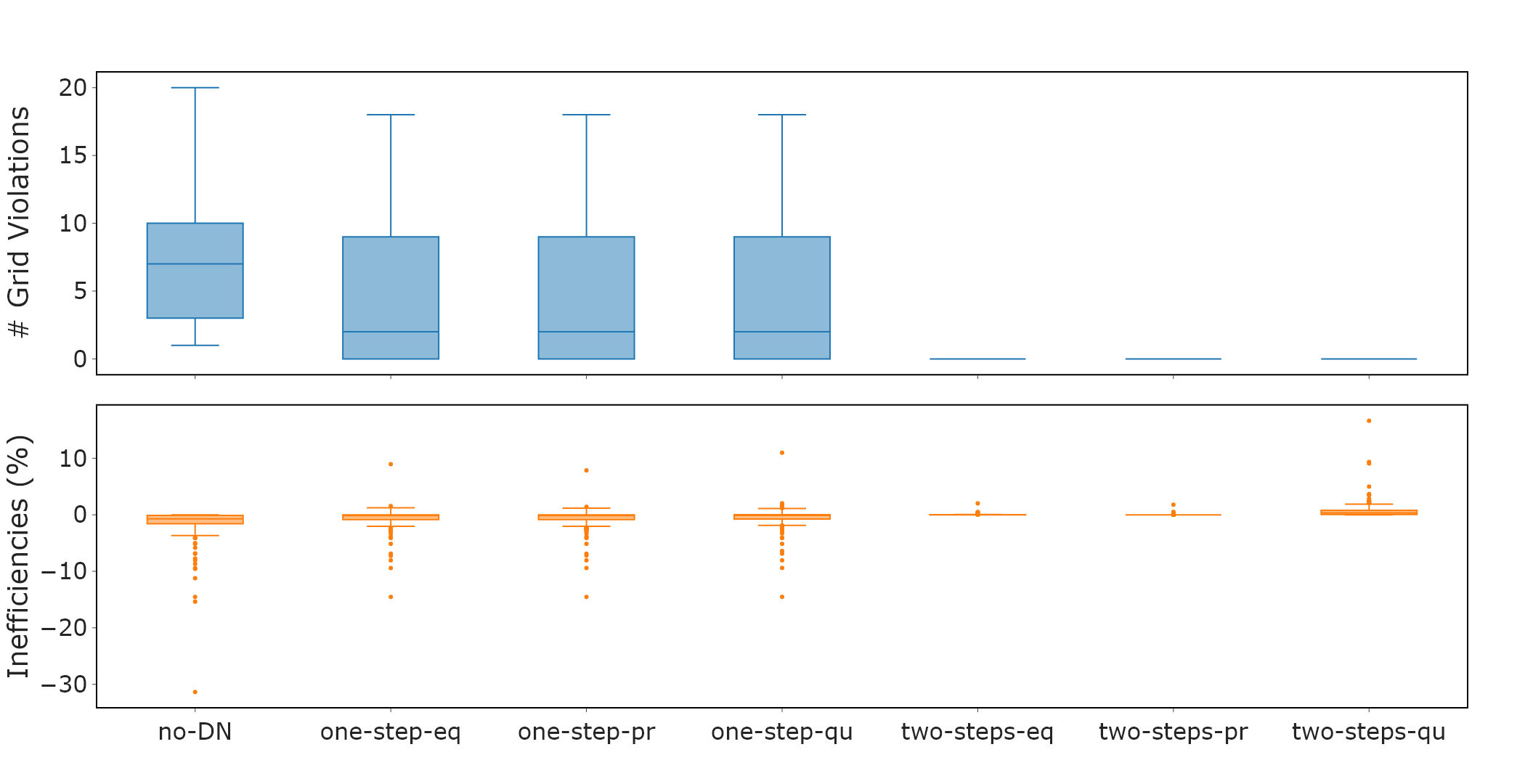}\vspace{-0.3cm}
\caption{(Top) Comparison of total numbers of violations for case set 1. (Bottom) Comparison of market inefficiencies for case set 1. Note that $\eta=0$ indicates the procurement cost is equal to that of the \emph{full-DN}.}
\label{fig:n_violations_inefficiencies}\vspace{-0.4cm}
\end{figure} 
The simulation results of case set 1 are shown in Figures~\ref{fig:n_violations_inefficiencies}\footnote{The \emph{full-DN} is not shown in the top plot as all results are grid-safe, neither in the bottom plot as all results equal 0 when applying equation \eqref{eq:inefficiency}.} and \ref{fig:unqualified_flex_case1}. From the top plot of Figure~\ref{fig:n_violations_inefficiencies}, we observe that Algorithm~\ref{alg:oe_liu2022} remarkably can ensure that the OE-based market-clearing outcomes are grid-safe, thus performing as well as the \emph{full-DN}. Unfortunately, this is not the case for Algorithm~\ref{alg:oe_petrou2021} for which the cleared flexibility can still cause grid violations. In the bottom plot of Figure~\ref{fig:n_violations_inefficiencies}, we can then observe that Algorithm~\ref{alg:oe_liu2022} results in non-negative inefficiencies, implying that the cleared bids are, albeit feasible, suboptimal. On the other hand, Algorithm~\ref{alg:oe_petrou2021} obtains lower procurement costs as a result of (partially) clearing unsafe flexibility. As such, there is a trade-off between procurement cost and grid-safety within the methods studied. Finally, Figure~\ref{fig:unqualified_flex_case1} shows the amount of flexibility that is disregarded due to the OE limits. Algorithm~\ref{alg:oe_liu2022} provides stricter limits for the downward flexibility than Algorithm~\ref{alg:oe_petrou2021}, which explains grid-safety and efficiency results. It is worth mentioning that, for Algorithm~\ref{alg:oe_liu2022}, despite the amount of unqualified flexibility being relatively high (in the range of $20\%$ on average), the loss in the market efficiency is not significant, which shows the effectiveness of the OE limits in restricting the distribution-level resources. We note that on average about $50\%$ of the flexibility needs of the TSO are satisfied by the DN resources. We also note that, in case set 1, almost no restriction is seen on the upward flexibility when applying the OE-methods (maximum of 0.06\% of unqualified upward bids is seen in Figure~\ref{fig:unqualified_flex_case1}). This is due to the case involving only shiftable loads. 

Next, we evaluate the sensitivity of Algorithms~\ref{alg:oe_liu2022}
and \ref{alg:oe_petrou2021} with respect to the weights. For Algorithm~\ref{alg:oe_liu2022}, the market efficiency is the best when price-based or equal-based weights are chosen conforming the objective of this weight rule (see bottom plot in Figure~\ref{fig:n_violations_inefficiencies}). The quantity-based weight rule results in more unqualified flexibility (see Figure~\ref{fig:unqualified_flex_case1}) and this can explain the worse performance in market efficiency compared to the price-based weight as unnecessary amount of flexibility is discarded (more restricted operational limits) resulting in larger inefficiency. On the other hand, Algorithm~\ref{alg:oe_petrou2021} is not sensitive to the weights as its performances are relatively similar with different weight rules imposed. Similarly to Algorithm~\ref{alg:oe_liu2022}, the quantity-based weight rule results in slightly larger unqualified flexibility. However, this does not impact the performance of Algorithm~\ref{alg:oe_petrou2021} on the other metrics.

\begin{figure}
\centering
\includegraphics[width=0.95\columnwidth]{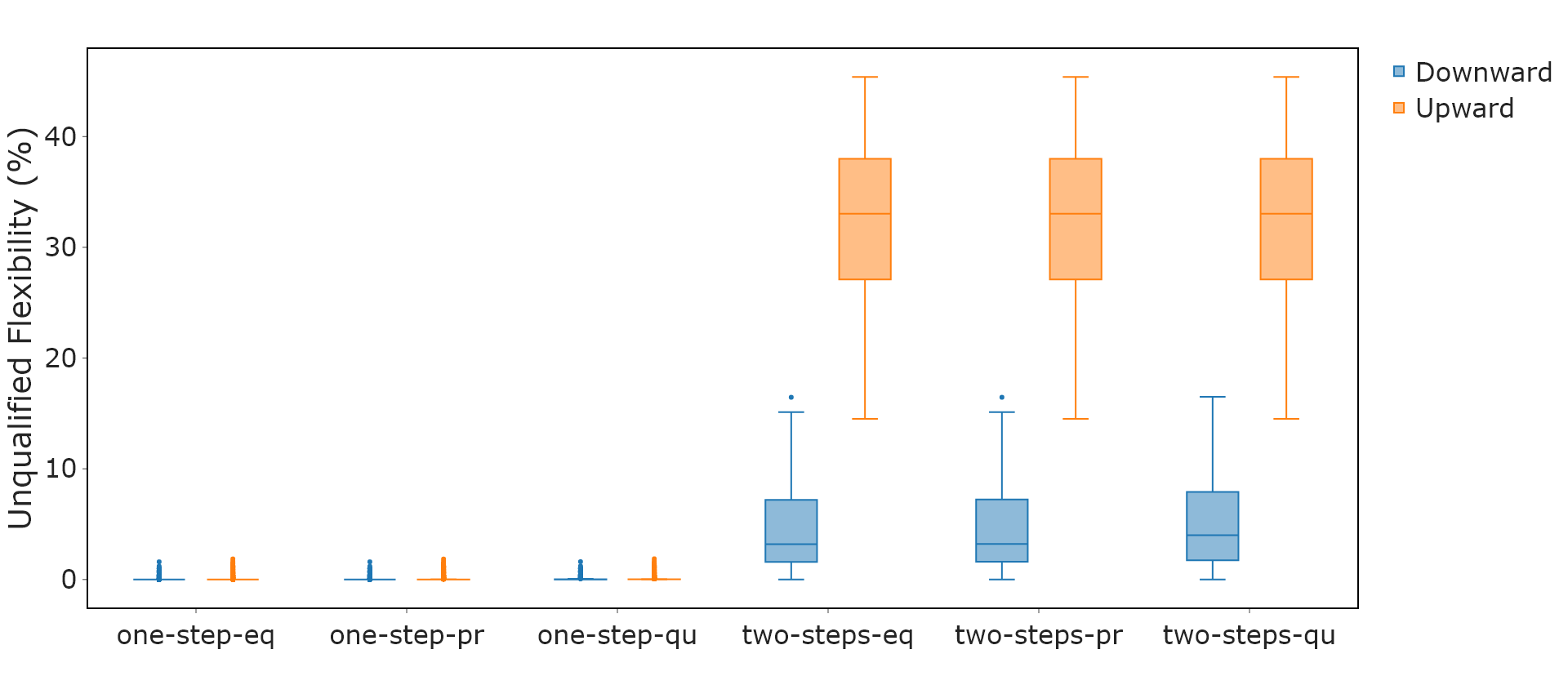}
\vspace{-0.3cm}
\caption{Comparison of unqualified flexibility for case set 2.}
\label{fig:unqualified_flex_case2}\vspace{-0.4cm}
\end{figure}

The results of case set 2 are in alignment with those of case set 1 for the number of grid violations and market inefficiencies. Therefore, we do not include the resulting plots due to space limitation. The only differences are in terms of the scale of results: 1) the number of grid violations can reach up to 35 in the \emph{no-DN} and Algorithm~\ref{alg:oe_petrou2021} scenarios; 2) the average market inefficiencies of Algorithm~\ref{alg:oe_liu2022} are 0, meaning that the OE-based market obtains cleared bids that are grid-safe and as efficient as the \emph{full-DN} market. This reinforces the efficiency of the two-step approach in Algorithm~\ref{alg:oe_liu2022}. Regarding unqualified bids, we notice that Algorithm~\ref{alg:oe_liu2022} restricts $30\%$, on average, of the upward flexibility (see Figure~\ref{fig:unqualified_flex_case2}), but also $5\%$, on average, of the downward flexibility, which is explained by the high amount of distributed resources available in this case. Although high, these limitations did not impact the market efficiency (as previously analyzed) implying that most of the unqualified flexibility would not have been cleared. Finally, we observe that, in most instances, Algorithm~\ref{alg:oe_petrou2021} does not impose restrictions to the resources as the OE limits are equal to their technical limits, resulting in similar performance as the \emph{no-DN} market in terms of grid violations and market efficiency.
\section{Conclusion} \label{sec:conclusion}
{Operating envelopes can be integrated into the flexibility market to act as a proxy of distribution network constraints with the aim of maintaining a safe operation of distribution systems. From our study, the two-step operating envelope calculation method, where the upward and downward operating envelopes are calculated separately outperforms the one-step method, where all operating envelopes are calculated simultaneously, as the former can produce grid-safe cleared flexibility at the cost of relatively low inefficiency while the latter can still result in grid violations. For the two-step algorithm, the price-based weight rule can be considered for a better performance in terms of market efficiency.
As our ongoing work, we are currently performing rigorous analysis to obtain a theoretical grid-safety guarantee of the two-step method and extending our numerical simulations by considering a nonlinear power flow model for grid-safety evaluation.
}
\vspace*{-4pt}
\bibliography{ref}

\begin{thebibliography}{10}

\bibitem{JRC130519}
D.~Koolen, M.~De~Felice, and S.~Busch, ``Flexibility requirements and the role
  of storage in future {European} power systems,'' tech. rep., Publications
  Office of the European Union, Luxembourg, 2023.

\bibitem{EUDSOcommissionregulation}
{Commission Regulation EU}, ``{EUDSO} {Entity} and {ENTSO-E} {DRAFT} {Proposal}
  for a {Network} {Code} on {Demand} {Response}.''
  \url{https://consultations.entsoe.eu/markets/public-consultation-networkcode-demand-response}.
\newblock Accessed on 11/04/2024.

\bibitem{marques2023grid}
L.~Marques~et al., ``Grid impact aware {TSO-DSO} market models for flexibility
  procurement: Coordination, pricing efficiency, and information sharing,''
  {\em IEEE Trans. on Power Systems}, vol.~38, no.~2, pp.~1920--1933, 2023.

\bibitem{petrou2021ensuring}
K.~Petrou~et al., ``Ensuring distribution network integrity using dynamic
  operating limits for prosumers,'' {\em IEEE Trans. on Smart Grid}, vol.~12,
  no.~5, pp.~3877--3888, 2021.

\bibitem{liu2021grid}
M.~Z. Liu~et al., ``Grid and market services from the edge: Using operating
  envelopes to unlock network-aware bottom-up flexibility,'' {\em IEEE Power
  and Energy Magazine}, vol.~19, no.~4, pp.~52--62, 2021.

\bibitem{liu2022using}
M.~Z. Liu~et al., ``Using {OPF}-based operating envelopes to facilitate
  residential {DER} services,'' {\em IEEE Trans. on Smart Grid}, vol.~13,
  no.~6, pp.~4494--4504, 2022.

\bibitem{8810638}
S.~Riaz and P.~Mancarella, ``On feasibility and flexibility operating regions
  of virtual power plants and {TSO/DSO} interfaces,'' in {\em 2019 IEEE Milan
  PowerTech}, pp.~1--6, 2019.

\bibitem{azim2024dynamic}
M.~I. Azim~et al., ``Dynamic operating envelope-enabled {P2P} trading to
  maximize financial returns of prosumers,'' {\em IEEE Trans. on Smart Grid},
  vol.~15, no.~2, pp.~1978--1990, 2024.

\bibitem{hoque2024dynamic}
M.~M. Hoque~et al., ``Dynamic operating envelope-based local energy market for
  prosumers with electric vehicles,'' {\em IEEE Trans. on Smart Grid}, vol.~15,
  no.~2, pp.~1712--1724, 2024.

\bibitem{sanjab2023joint}
A.~Sanjab, L.~Marques, H.~Gerard, and K.~Kessels, ``Joint and sequential
  {DSO-TSO} flexibility markets: efficiency drivers and key challenges,'' in
  {\em CIRED}, vol.~2023, pp.~3138--3143, IET, 2023.

\bibitem{marques2023flexibility}
L.~Marques, A.~Sanjab, and T.~Cuypers, ``Flexibility service providers' gaming
  potential and its impact on {TSO-DSO} coordinated markets,'' in {\em SEST
  Conference}, pp.~1--6, IEEE, 2023.

\bibitem{sanjab2021linear}
A.~Sanjab, Y.~Mou, A.~Virag, and K.~Kessels, ``A linear model for distributed
  flexibility markets and {DLMPs}: A comparison with the {SOCP} formulation,''
  in {\em CIRED}, vol.~2021, pp.~3181--3185, 2021.

\end{thebibliography}
\bibliographystyle{ieeetr}
\end{document}